\documentclass{amsart}

\newtheorem{theorem}{Theorem}
\newcommand{\bt}{\begin{theorem}}
\newcommand{\et}{\end{theorem}}
\newtheorem{lemma}{Lemma}
\newcommand{\bl}{\begin{lemma}}
\newcommand{\el}{\end{lemma}}
\newtheorem{corollary}{Corollary}
\newcommand{\bc}{\begin{corollary}}
\newcommand{\ec}{\end{corollary}}
\newtheorem{problem}{Problem}
\newcommand{\bprob}{\begin{problem}}
\newcommand{\eprob}{\end{problem}}

\newcommand{\beq}{\begin{equation}}
\newcommand{\eeq}{\end{equation}}
\newcommand{\benum}{\begin{enumerate}}
\newcommand{\eenum}{\end{enumerate}}
\newcommand{\N}{\ensuremath{ \mathbf N }}
\newcommand{\Z}{\ensuremath{\mathbf Z}}

\newcommand{\R}{\ensuremath{\mathbf R}}

\newcommand{\mcy}{\ensuremath{ \mathcal Y}}

\newcommand{\mba}{\ensuremath{ \mathbf a}}
\newcommand{\mbb}{\ensuremath{ \mathbf b}}

\newcommand{\mbv}{\ensuremath{ \mathbf v}}

\newcommand{\mbx}{\ensuremath{ \mathbf x}}
\newcommand{\mby}{\ensuremath{ \mathbf y}}

\newcommand{\Rn}{\ensuremath{ \mathbf{R}^n }}

\newcommand{\bsmallmat}{\left(\begin{smallmatrix}}
\newcommand{\esmallmat}{\end{smallmatrix}\right)}

\newcommand{\bmat}{\left(\begin{matrix}}
\newcommand{\emat}{\end{matrix}\right)}

\DeclareMathOperator{\qqand}{\qquad\text{and}\qquad}

\DeclareMathOperator{\vectorsmallan}{\left( \begin{smallmatrix} a_1 \\ \vdots \\ a_n \end{smallmatrix}\right)}
\DeclareMathOperator{\vectorsmallbn}{\left( \begin{smallmatrix} b_1 \\ \vdots \\ b_n \end{smallmatrix}\right)}

\DeclareMathOperator{\vectorsmallxn}{\left( \begin{smallmatrix} x_1 \\ \vdots \\ x_n \end{smallmatrix}\right)}

\usepackage{amssymb,latexsym}
\usepackage[all]{xy}

\title{Underapproximation by Egyptian fractions}
\author{Melvyn B. Nathanson}
\address{Department of Mathematics\\Lehman College (CUNY)\\Bronx, NY 10468} 
\email{melvyn.nathanson@lehman.cuny.edu} 
\date{\today}



\thanks{Supported in part by a grant from the PSC-CUNY Research Award Program.}

\begin{document}
\maketitle

\begin{abstract}
An increasing sequence $(x_i)_{i=1}^n$ of positive integers is an  $n$-term Egyptian 
underapproximation sequence of $\theta \in (0,1]$ if $\sum_{i=1}^n \frac{1}{x_i}  < \theta$.
A greedy algorithm constructs an $n$-term underapproximation sequence of $\theta$.  
For some but not all  numbers $\theta$, the greedy algorithm gives a unique best $n$-term 
underapproximation sequence for all $n$.  An infinite set of rational numbers is constructed for which 
the greedy underapproximations are best, and numbers for which the greedy algorithm 
is not best are also studied.     
\end{abstract}

\section{The greedy underapproximation algorithm}

An \emph{Egyptian fraction}\index{Egyptian fraction} is a fraction of the form $1/x$, 
where $x$ is a positive integer.   
Let $\theta \in (0,1]$. 
A finite sequence $(x_i)_{i=1}^n$ of integers is an \emph{$n$-term Egyptian 
underapproximation sequence}\index{ Egyptian underapproximation} 
of $\theta$ if 
\[
2 \leq x_1 \leq x_2 \leq \cdots \leq x_n 
\]
and
\[
\sum_{i=1}^n \frac{1}{x_i}  < \theta.
\]  
For example, $(2,3,7,43)$ is a 4-term underapproximation sequence of 1.  
If $x$ is an integer such that $x > n/\theta$, then 
\[
\sum_{i=1}^n \frac{1}{x+i} < \theta
\] 
and $(x+i)_{i=1}^n$ is  is an $n$-term Egyptian underapproximation sequence of $\theta$.  

An infinite sequence $(x_i)_{i=1}^{\infty}$ of integers is an \emph{infinite Egyptian 
underapproximation sequence} of $\theta$ if 
the finite sequence $(x_i)_{i=1}^n$  is an  $n$-term Egyptian 
underapproximation sequence of $\theta$ for all $n \geq 1$.

For all $\theta \in (0,1]$, there is a  unique positive integer  $G(\theta) = a$ such that 
\[
a \geq 2       \qqand \frac{1}{a} < \theta \leq \frac{1}{a-1}.
\]
Thus, $G(\theta)$ is the smallest positive integer such that  Egyptian fraction $1/G(\theta)$  
underapproximates $\theta$.  
Equivalently, 
\[
a \leq \frac{1}{\theta} + 1 < a+1 
\]
and so\footnote{ The \emph{greatest integer function}\index{greatest integer function}  of the real 
number $w$, also called the \emph{floor} of $w$, is the unique integer $\ell$ 
such that $\ell \leq w < \ell+1$.  
We write $\lfloor w \rfloor = \ell$.  
The \emph{ceiling} of $w$, denoted $\lceil w \rceil$, is the unique integer $m$ 
such that $m \geq w > m-1$.  Define the interval $(t_1, t_2] = \{ t \in \R: t_1 < t \leq t_2 \}$.}
\[
G(\theta) = \left\lfloor  \frac{1}{\theta} \right\rfloor + 1.
\]

For all $\theta \in (0,1]$, the 
\emph{greedy underapproximation algorithm}\index{greedy underapproximation algorithm}
\index{underapproximation algorithm} applied to $\theta$ constructs an infinite sequence of integers 
$(a_i)_{i=1}^{\infty}$ as follows:
\beq                    \label{Egyptian:underapproxIneqality-0} 
a_1 = G(\theta) \geq 2 
\eeq
and, for all  $i \geq 1$ and integers $a_1,a_2,\ldots, a_i$,
\[
a_{i+1} = G\left( \theta - \sum_{j =1}^i \frac{1}{a_j} \right). 
\]
Thus,  
\beq                    \label{Egyptian:underapproxIneqality-1} 
\frac{1}{a_{i+1}} <  \theta - \sum_{j=1}^i \frac{1}{a_j} \leq \frac{1}{a_{i+1} -1}. 
\eeq
Equivalently, 
\[
\sum_{j =1}^{i+1} \frac{1}{a_j} <  \theta \leq  \sum_{j =1}^{i} \frac{1}{a_j} + \frac{1}{a_{i+1} - 1}.
\]
We call $(a_i)_{i=1}^{\infty}$ the 
\emph{infinite greedy underapproximation sequence}\index{greedy underapproximation sequence!infinite} 
of $\theta$ and $(a_i)_{i=1}^n$ the 
\emph{$n$-term greedy underapproximation sequence}\index{greedy underapproximation sequence! $n$-term} 
of $\theta$.  The rational number $\sum_{i=1}^n 1/a_i$ is the 
\emph{$n$-term greedy underapproximation} of $\theta$.

 Let $ i \geq 1$.  
Inequality~\eqref{Egyptian:underapproxIneqality-1} implies that 
\begin{align*}
\frac{1}{a_{i+1}} & <  \theta - \sum_{j=1}^i \frac{1}{a_j}  
= \left( \theta - \sum_{j=1}^{i-1} \frac{1}{a_j}  \right) - \frac{1}{a_{i} } \\
& \leq \frac{1}{a_{i} -1}  - \frac{1}{a_{i} } 
 = \frac{1}{a_{i} (a_{i} -1)}
\end{align*}
and so 
\beq                    \label{Egyptian:underapproxIneqality-2} 
a_{i+1} \geq a_i^2 - a_i + 1. 
\eeq
It follows from~\eqref{Egyptian:underapproxIneqality-1} 
and~\eqref{Egyptian:underapproxIneqality-2} that $(a_i)_{i=1}^{\infty}$, 
the infinite greedy underapproximation sequence of $\theta$, is a 
strictly increasing sequence of positive integers and  that 
\[
\sum_{i=1}^{\infty} \frac{1}{a_i} = \theta.  
\]

Here is a classical example of Egyptian underapproximation.  
\emph{Sylvester's sequence}~\cite{sylv80}   
is the sequence of positive integers $(s_i)_{i=1}^{\infty}$ 
constructed recursively by the following rule: 
\beq             \label{Egyptian:Sylvester} 
s_1 = 2 \qqand 
s_{i+1} =   \prod_{j=1}^i s_j   \ + 1 
\eeq 
for all $i \geq 1$.  We have
\begin{align*}
s_1 & = 2 \\
s_2 & =  3 \\
s_3 & =  7 \\
s_4 & =  43 \\
s_5 & =  1807 \\
s_6 & =  3263443 \\
s_7 & =  10650056950807 \\
s_8 & =   113423713055421844361000443 \\
s_9 & = 12864938683278671740537145998360961546653259485195807.  
\end{align*}
Sylvester's sequence is sequence A000058 in the OEIS.
By Corollary~\ref{Egyptian:corollary:Sylvester},  
Sylvester's sequence $(s_i)_{i=1}^{\infty}$ 
is the infinite greedy underapproximation sequence of $\theta = 1$.

The following theorem constructs a  set of rational numbers whose infinite 
greedy approximation sequences generalize Sylvester's sequence.  

\bt                \label{Egyptian:theorem:pq-greedy-sequence}
Let $\theta = p/q \in (0,1]$, where $p$ and $q$ are positive integers such that $p$ divides $q+1$, 
and let $(a_i)_{i=1}^{\infty}$ be the infinite greedy underapproximation 
sequence of $\theta$.   Then 
\[
a_1 = \frac{q+1}{p}
\]
and, for all $k \geq 1$,  
\[
a_{k+1} = q\prod_{i=1}^k a_i + 1 
\]
and 
\[
\frac{p}{q}  = \sum_{i=1}^k \frac{1}{a_i} + \frac{1}{q\prod_{i=1}^k a_i}. 
\] 
\et

\begin{proof}
The proof is by induction on $k$.  
Let $q+1 = pt$.  We have 
\[
\frac{1}{t}  =  \frac{p}{q+1} <  \frac{p}{q}  = \frac{p}{pt-1} \leq \frac{1}{t-1}
\]
and so 
\[
a_1 = G\left( \frac{p}{q} \right) = t = \frac{q+1}{p}.
\]
It follows that 
\[
\frac{p}{q} - \frac{1}{a_1} = \frac{p}{q} - \frac{1}{t} =  \frac{1}{qt} = \frac{1}{qa_1} 
\]
and  so 
\[
a_2 =  G\left( \frac{p}{q} - \frac{1}{a_1} \right) = qa_1 + 1.  
\]
We obtain  
\[
\frac{p}{q} - \frac{1}{a_1} - \frac{1}{a_2} = \frac{1}{qa_1}  - \frac{1}{qa_1 + 1}  
=  \frac{1}{q a_1 (qa_1 + 1)}  = \frac{1}{q a_1 a_2 } 
\] 
and  so 
\[
a_3 =  G\left( \frac{p}{q} - \frac{1}{a_1}  - \frac{1}{a_2} \right) = qa_1 a_2 + 1.  
\]
Let $k \geq 2$.  If  
\[
\frac{p}{q} - \sum_{i=1}^k \frac{1}{a_i} = \frac{1}{q\prod_{i=1}^k a_i} 
\]
then 
\[
a_{k+1} =  G\left( \frac{p}{q} -  \sum_{i=1}^k \frac{1}{a_i}\right) = q\prod_{i=1}^k a_i + 1
\]
and 
\begin{align*}
\frac{p}{q} - \sum_{i=1}^{k+1} \frac{1}{a_i} 
& = \frac{p}{q} - \sum_{i=1}^{k} \frac{1}{a_i} - \frac{1}{a_{k+1}} = \frac{1}{q\prod_{i=1}^k a_i} - \frac{1}{ q\prod_{i=1}^k a_i + 1} \\ 
& = \frac{1}{q\prod_{i=1}^k a_i \left(  q\prod_{i=1}^k a_i + 1\right) } \\
& = \frac{1}{ q \prod_{i=1}^{k+1} a_i}. 
\end{align*}
This completes the proof.  
\end{proof} 

\bc                  \label{Egyptian:corollary:Sylvester}
Sylvester's sequence is the infinite greedy underapproximation sequence for $\theta = 1$.
\ec

\section{A criterion for greedy underapproximation}  \label{Egyptian:section:criterion}

\bt                                 \label{Egyptian:theorem:greedy-n-term-condition} 
Let $(a_i)_{i=1}^n$ be a sequence of integers such that 
\[
a_1 \geq 2 \qqand a_{i+1} \geq a_i^2 - a_i + 1
\]
for all $i  = 1,\ldots, n-1$.  The sequence 
$(a_i)_{i=1}^n$ is the $n$-term greedy underapproximation sequence of the real number $\theta$ 
if and only if 
\beq                    \label{Egyptian:greedy-n-term-interval} 
\theta \in \left( \sum_{i=1}^n \frac{1}{a_i} , \  \sum_{i=1}^{n-1} \frac{1}{a_i} + \frac{1}{a_n -1} \right].
\eeq
\et

\begin{proof}
If $(a_i)_{i=1}^n$ is the $n$-term greedy underapproximation sequence of  $\theta$, then 
\[
\frac{1}{a_n} < \theta - \sum_{i=1}^{n-1} \frac{1}{a_i} \leq  \frac{1}{a_n-1}
 \]
and so $\theta$  is in the interval~\eqref{Egyptian:greedy-n-term-interval}. 

To prove the converse, we observe that, for all $i = 1,\ldots, n-1$, 
the inequality $a_{i+1} \geq a_i^2 - a_i + 1$ implies that 
\[
\frac{1}{a_i} + \frac{1}{a_{i+1}-1} \leq  \frac{1}{a_i-1}. 
\]
It follows that, for all  $k = 1, \ldots, n$, we have 
\begin{align*}
\sum_{i=1}^{n-1} \frac{1}{a_i} + \frac{1}{a_n -1}
& = \sum_{i=1}^{n-2} \frac{1}{a_i} + \frac{1}{a_{n -1}} + \frac{1}{a_n -1} 
 \leq \sum_{i=1}^{n-2} \frac{1}{a_i} + \frac{1}{a_{n -1} -1}  \\
& \leq  \cdots  \leq \sum_{i=1}^{k-1} \frac{1}{a_i} + \frac{1}{a_k-1} 
\end{align*} 
and so 
\[
\sum_{i=1}^{k} \frac{1}{a_i} \leq  \sum_{i=1}^{n} \frac{1}{a_i} 
<  \sum_{i=1}^{n-1} \frac{1}{a_i} + \frac{1}{a_n -1} 
\leq \sum_{i=1}^{k-1} \frac{1}{a_i} + \frac{1}{a_k-1}
\]
If $\theta$  is in the interval~\eqref{Egyptian:greedy-n-term-interval}, 
then for all  $k = 1, \ldots, n$ we have  
\[
\sum_{i=1}^{k} \frac{1}{a_i} 
<  \theta 
\leq \sum_{i=1}^{k-1} \frac{1}{a_i} + \frac{1}{a_k-1}. 
\]
Equivalently,  
\[
\frac{1}{a_k} < \theta - \sum_{i=1}^{k-1} \frac{1}{a_i} \leq  \frac{1}{a_k-1}
 \]
and 
\[
a_k  = G\left( \theta - \sum_{i=1}^{k-1} \frac{1}{a_i}  \right).
\]
Thus, $(a_i)_{i=1}^n$ is the $n$-term greedy underapproximation sequence of  $\theta$. 
This completes the proof.  
\end{proof}

\bc                   \label{Egyptian:corollary:GUAsequence-1}
Let $\theta \in (0,1]$. 
The pair of integers $(a_1,a_2)$ with $2 \leq a_1 \leq a_2$  
is the 2-term greedy underapproximation sequence of $\theta$ 
if and only if  $a_2 \geq a_1^2 - a_1 + 1$ and 
\[
\frac{1}{a_1} + \frac{1}{a_2} < \theta \leq \frac{1}{a_1} + \frac{1}{a_2-1}.
\]
\ec

\bc                                     \label{Egyptian:corollary:TheConverseAlgorithm}  
Let $(a_i)_{i=1}^{\infty}$ be a sequence of integers  such that 
\[
a_1 \geq 2 \qqand a_{i+1} \geq a_i^2 - a_i + 1
\]
for all $i \geq 1$.
The infinite series 
\[
\sum_{i=1}^{\infty} \frac{1}{a_i} 
\]
converges to a number $\theta \in (0,1]$, and $(a_i)_{i=1}^{\infty}$ is the infinite greedy 
underapproximation sequence of $\theta$.   
\ec

\section{Best Egyptian approximation}            \label{Egyptian:section:BestApproximation}

Let $E_n$ be the set of all sequences $(x_i)_{i=1}^n$ of integers such that 
\[
  2 \leq x_1 \leq x_2 \leq \cdots \leq x_n.  
\]
For $\theta \in (0,1]$, let $U_n(\theta)$ be the set of all $n$-term Egyptian 
underapproximation sequences of $\theta$.  Thus,
\[
U_n(\theta) = \left\{
(x_i)_{i=1}^n \in E_n: \sum_{i=1}^n \frac{1}{x_i} < \theta 
\right\}.
\]
Let
\[
u_n(\theta) = \sup \left\{  \sum_{i=1}^n \frac{1}{x_i} : (x_i)_{i=1}^n \in U_n(\theta) \right\}. 
\]
We call $u_n(\theta)$ the \emph{best $n$-term Egyptian underapproximation of $\theta$}.  

If $(x_i)_{i=1}^n \in U_n(\theta)$, then $\sum_{i=1}^n 1/x_i < \theta$ and so 
$u_n(\theta) \leq \theta$.  
We shall prove (Theorem~\ref{Egyptian:theorem:BestUnderapproximation}) 
that there is a sequence $(b_i)_{i=1}^n \in U_n(\theta)$ 
such that $u_n(\theta) = \sum_{i=1}^n 1/b_i$ and so 
$u_n(\theta)$ is a rational number that is strictly less than $\theta$. 
We shall also construct examples to prove that the $n$-term greedy underapproximation 
of $\theta$ is not necessarily the best $n$-term Egyptian underapproximation 
and that there is not necessarily a unique sequence 
that is the best $n$-term Egyptian underapproximation of $\theta$.

\bt              \label{Egyptian:theorem:BestUnderapproximation} 
Let $\theta \in (0,1]$.  For all $n \geq 1$, there is a sequence $ (b_i)_{i=1}^n \in U_n(\theta)$ such that 
\[
u_n(\theta) =  \sum_{i=1}^n \frac{1}{b_i} < \theta. 
\]
Thus, the best $n$-term underapproximation $u_n(\theta)$ is rational.
\et

\begin{proof} 
For $n=1$, we have 
\[
U_1(\theta) = \left\{(x_1) : x_1 \geq a_1 = G(\theta) \right\}. 
\]
Setting $b_1 = a_1$ gives $u_1(\theta) = 1/a_1 = 1/b_1 < \theta$.  

Let $n \geq 2$.
Choose an $n$-tuple  $\left( c^{(1)}_i \right)_{i=1}^n \in U_n(\theta)$.   We have 
\[
  2 \leq c^{(1)}_1 \leq c^{(1)}_2 \leq \cdots \leq c^{(1)}_n \qqand \sum_{i=1}^n \frac{1}{c^{(1)}_i} < \theta. 
\]
If $\left( x_i \right)_{i=1}^n \in U_n(\theta)$ and 
\[
x_1 \geq nc^{(1)}_1 = x_1^*
\]
then the inequality  $x_1 \leq x_2 \leq \cdots \leq x_n$ 
implies that 
\[
\sum_{i=1}^n \frac{1}{x_i} \leq \frac{n}{x_1} \leq \frac{1}{c^{(1)}_1} < \sum_{i=1}^n \frac{1}{c^{(1)}_i} < \theta.
\]  
Thus, $\left( c^{(1)}_i \right)_{i=1}^n$ is a larger $n$-term Egyptian underapproximation 
of $\theta$ than $\left( x_i \right)_{i=1}^n$. 
Let 
\[
U_n^{(1)}(\theta) =  \left\{  \left( x_i \right)_{i=1}^n \in U_n(\theta) 
 \text{ and }  x_1 < x_1^*  \right\}. 
\]
We have   
\begin{align*}
u_n(\theta) 
& = \sup \left\{  \sum_{i=1}^n \frac{1}{x_i}    : (x_i)_{i=1}^n \in U_n(\theta) \right\} \\
& = \sup \left\{  \sum_{i=1}^n \frac{1}{x_i}   :  (x_i)_{i=1}^n \in U_n(\theta) 
 \text{ and }  x_1 < x_1^* \right\} \\
 & = \sup \left\{  \sum_{i=1}^n \frac{1}{x_i}   :  (x_i)_{i=1}^n \in U_n^{(1)}(\theta) \right\}.
\end{align*}

Let $k \in \{1,\ldots, n-1\}$ and let $x_1^*,\ldots, x_k^*$ be positive integers  
 such that 
 \[
u_n(\theta) 
 = \sup \left\{  \sum_{i=1}^n \frac{1}{x_i}   :  \left( x_i \right)_{i=1}^n \in U_n(\theta) 
 \text{ and }  x_i < x_i^* \text{ for all } i=1,\ldots, k  \right\}.
 \]
 Let 
\[
U_n^{(k)}(\theta) =  \left\{  \left( x_i \right)_{i=1}^n \in U_n(\theta) :   x_i < x_i^* \text{ for all } i=1,\ldots, k    \right\}. 
\]
Thus, 
\[
u_n(\theta) 
 = \sup \left\{  \sum_{i=1}^n \frac{1}{x_i}   :  (x_i)_{i=1}^n \in U_n^{(k)}(\theta) \right\}.
 \]

Let $\mcy(k,n)$ be the finite set of all $k$-tuples of positive integers  
$\mby = \left( y_i \right)_{i=1}^k$ such that 
\benum
\item[(i)]
$y_i < x_i^*$ for all $i=1,\ldots, k$,  and 
\item[(ii)]  
there exists an $n$-tuple $\left( x_i \right)_{i=1}^n \in U_n^{(k)}(\theta)$ such that 
$x_i = y_i$ for all $i = 1,\ldots, k$. 
\eenum
For each $k$-tuple $\mby = \left( y_i \right)_{i=1}^k \in \mcy(k,n)$, 
let $U_n^{ (\mby)}(\theta)$ be the nonempty set of all  $n$-tuples  
$\left( x_i \right)_{i=1}^n \in U_n^{(k)}(\theta)$ 
such that $x_i = y_i$ for all $i = 1,\ldots, k$.  We have 
\[
U_n^{(k)}(\theta) = \bigcup_{\mby \in \mcy(k,n)} U_n^{ (\mby)}(\theta). 
\]

For all $\mby \in \mcy(k,n)$, choose  an $n$-tuple 
$\left(c^{(\mby)}_i \right)_{i=1}^n \in U_n^{(\mby)}(\theta)$.  
If $\left( x_i \right)_{i=1}^n \in U_n^{(\mby)}(\theta)$ and   
\[
x_{k+1} \geq (n-k)c_{k+1}^{(\mby)}
\]
then $x_i = y_i = c^{(\mby)}_i $ for all $i=1,\ldots, k$ and 
\begin{align*} 
\sum_{i=1}^n \frac{1}{x_i} & = \sum_{i=1}^k \frac{1}{ c^{(\mby)}_i }  + \sum_{i=k+1}^n \frac{1}{x_i} \\
& \leq \sum_{i=1}^k \frac{1}{ c^{(\mby)}_i }  + \frac{n-k}{x_{k+1}} \\
& \leq \sum_{i=1}^{k+1} \frac{1}{ c^{(\mby)}_i }  \leq \sum_{i=1}^{n} \frac{1}{ c^{(\mby)}_i } \\
& < \theta
\end{align*} 
and so the $n$-term  Egyptian underapproximation of $\theta$ by  $\left( x_i \right)_{i=1}^n$ 
is no larger than the $n$-term  Egyptian underapproximation of $\theta$ 
by  $\left(c^{(\mby)}_i \right)_{i=1}^n$. 
Therefore, 
\begin{align*}
\sup  & \left\{  \sum_{i=1}^n \frac{1}{x_i}  :  (x_i)_{i=1}^n \in U_n^{(\mby)}(\theta) \right\} \\
& =  \sup \left\{  \sum_{i=1}^n \frac{1}{x_i}  :  (x_i)_{i=1}^n \in U_n^{(\mby)}(\theta) 
\text{ and } x_{k+1} < (n-k)c_{k+1}^{(\mby)}\right\}. 
\end{align*}
Let 
\[
x_{k+1}^* = \max\left\{(n-k)c_{k+1}^{(\mby)} : \mby \in \mcy(k,n) \right\}.  
\]
It follows that 
\begin{align*}
\sup  & \left\{  \sum_{i=1}^n \frac{1}{x_i}  :  (x_i)_{i=1}^n \in U_n^{(k)}(\theta) \right\} \\
& =  \sup \left\{  \sum_{i=1}^n \frac{1}{x_i}  :  (x_i)_{i=1}^n \in U_n^{(k)}(\theta) 
\text{ and } x_i < x_i^* \text{ for all } i = 1,\ldots, k+1\right\}. 
\end{align*}
Continuing inductively, we obtain positive integers $x_1^*,\ldots, x_n^*$ such that 
\begin{align*}
u_n(\theta) 
& = \sup \left\{  \sum_{i=1}^n \frac{1}{x_i} :  (x_i)_{i=1}^n \in U_n(\theta) \right\}  \\
&  = \sup \left\{  \sum_{i=1}^n \frac{1}{x_i} :  (x_i)_{i=1}^n \in U_n(\theta) 
  \text{ and }  x_i < x_i^* \text{ for all } i = 1,\ldots, n \right\} \\
  & =  \sup \left\{  \sum_{i=1}^n \frac{1}{x_i} :  (x_i)_{i=1}^n \in U_n^{(n)}(\theta) \right\} 
\end{align*}
where
\[
U_n^{(n)}(\theta) =  \left\{   (x_i)_{i=1}^n \in U_n(\theta) :   x_i < x_i^* \text{ for all } i = 1,\ldots, n \right\}. 
\]
The set $U_n^{(n)}(\theta)$ is finite and so 
there exists $ (b_i)_{i=1}^n \in U_n^{(n)}(\theta) \subseteq U_n(\theta)$ such that 
\[
u_n(\theta) =  \sum_{i=1}^n \frac{1}{b_i} < \theta.
\]
This completes the proof. 
\end{proof}

\section{When greedy  is best}              \label{Egyptian:section:bb}

It had been conjectured by Miller~\cite{mill19} and   Kellogg~\cite{kell21}
and then proved by Curtiss~\cite{curt22} and Takenouchi~\cite{take21}  
that, for every  positive integer $n$, the $n$-tuple of Sylvester numbers 
$(s_i)_{i=1}^n$  is the unique best $n$-term Egyptian fraction underapproximation of 1.  Equivalently, 
if $(x_1,\ldots, x_n) \in U_n(1)$ and 
\[
\sum_{i=1}^n \frac{1}{s_i} \leq \sum_{i=1}^n \frac{1}{x_i} < 1
\]
then $x_i = s_i$ for all $i=1,\ldots, n$.  There is also a recent proof by Soundararajan~\cite{soun05}.

In this section we generalize this result.  
We construct an infinite set of rational numbers 
whose infinite greedy underapproximation sequences can be expicitly computed, and for which, for every $n$, 
the $n$-term greedy underapproximation sequence is the unique best $n$-term underapproximation 
by Egyptian fractions.  We use the method of Soundararajan~\cite{soun05}, which is based 
on the following inequality.  

\bt                                                          \label{Egyptian:theorem:MuirheadCorollary}                     
If  $(x_i)_{i=m+1}^n$ and $(a_i)_{i=m+1}^n$ are   increasing sequences of positive numbers 
such that  $(x_i)_{i=m+1}^n \neq (a_i)_{i=m+1}^n$ and 
\[
\prod_{i=m+1}^{m+k} a_i \leq \prod_{i=m+1}^{m+k} x_i 
\]
for all $k = 1,\ldots, n-m$, then 
\[
\sum_{i= m+1}^n \frac{1}{x_i} < \sum_{i= m+1}^n  \frac{1}{a_i} .
\]
\et

\begin{proof}
This inequality is a corollary of Muirhead's inequality (see Nathanson~\cite{nath22}).  
A nice direct proof  due to Ambro and Barc\u{a}u~\cite{ambr-barc15} is given in the 
Appendix.
\end{proof}

\bt                          \label{Egyptian:theorem:pq-greedy}
Let $\theta = p/q \in (0,1]$, where $p$ and $q$ are positive integers such that $p$ divides $q+1$, 
and let $(a_i)_{i=1}^{\infty}$ be the infinite greedy underapproximation sequence of $\theta$.   
For every positive integer $n$, if $(x_i)_{i=1}^n$ is an $n$-term Egyptian 
underapproximation  sequence of $\theta$  such that 
\beq                              \label{Egyptian:pq-greedy}
\sum_{i=1}^n \frac{1}{a_i} \leq \sum_{i=1}^n \frac{1}{x_i} < \frac{p}{q} 
\eeq
then $x_i = a_i$ for all $i = 1,\ldots, n$. 
\et

\begin{proof}
The proof is by induction on $n$.  For $n=1$, the greedy algorithm gives  
\[
\frac{1}{a_1} \leq \frac{1}{x_1}  < \theta \leq \frac{1}{a_1 -1} 
\]
and so $x_1 = a_1$.  Thus,  the Theorem is true for $n=1$. 

Let $n \geq 2$ and assume that the Theorem is true for all increasing sequences 
$(x_i)_{i=1}^{m}$ with $m < n$.  
Let $(x_i)_{i=1}^n$ be an increasing sequence that satisfies~\eqref{Egyptian:pq-greedy}.  
Inequality~\eqref{Egyptian:pq-greedy} and 
Theorem~\ref{Egyptian:theorem:pq-greedy-sequence} give  
\[
0 <  \frac{p}{q}  -  \sum_{i=1}^n \frac{1}{x_i} \leq \frac{p}{q}  - \sum_{i=1}^n \frac{1}{a_i}  =  \frac{1}{q\prod_{i=1}^n a_i}.  
\]
A common denominator of the $n+1$ fractions $p/q$, $1/x_1$, \ldots, $1/x_n$ is $q\prod_{i=1}^n x_i$, 
and so there is a positive integer $r$ such that 
\[
0 < \frac{1}{q\prod_{i=1}^n x_i} \leq \frac{r}{q\prod_{i=1}^n x_i} 
=  \frac{p}{q}  -  \sum_{i=1}^n \frac{1}{x_i} \leq   \frac{1}{q\prod_{i=1}^n a_i}.  
\]
This implies  
\[
\prod_{i=1}^n a_i \leq \prod_{i=1}^n x_i.
\]

Let $m$ be the largest integer $\leq n-1$ such that 
\beq                              \label{Egyptian:m+1-product}
\prod_{i=m+1}^n a_i \leq \prod_{i=m+1}^n x_i.
\eeq 
We shall prove that 
\beq                              \label{Egyptian:j-product}
\prod_{i=m+1}^{m+j} a_i \leq \prod_{i=m+1}^{m+j} x_i
\eeq
for all $j \in \{1,\ldots, n-m-1\}$.
If not, then there exists $k \in \{1,\ldots, n-m-1\}$ such that 
\[
\prod_{i=m+1}^{m+k} x_i <\prod_{i=m+1}^{m+k} a_i.  
\]
It follows from~\eqref{Egyptian:m+1-product} that   
\[
\prod_{i=m+k+1}^n a_i \leq \frac{ \prod_{i=m+1}^n x_i}{\prod_{i=m+1}^{m+k} a_i} 
=  \left( \frac{ \prod_{i=m+1}^{m+k} x_i}{\prod_{i=m+1}^{m+k} a_i} \right)  \prod_{i=m+k+1}^n x_i 
< \prod_{i=m+k+1}^n x_i
\]
which contradicts the maximality of $m$.  
This proves~\eqref{Egyptian:j-product}.

Suppose  that 
 $a_i \neq x_i$ for some $i \in \{ m+1,\ldots, n\}$.   
Applying Theorem~\ref{Egyptian:theorem:MuirheadCorollary} 
to the distinct  increasing sequences $(a_i)_{i=m+1}^n$ and $(x_i)_{i=m+1}^n$, we obtain 
\beq                              \label{Egyptian:m-product}
\sum_{i=m+1}^{n} \frac{1}{x_i} < \sum_{i=m+1}^{n} \frac{1}{a_i}. 
\eeq 
From inequality~\eqref{Egyptian:pq-greedy} we have $1 \leq m \leq n-1$, and so 
\begin{align*} 
\sum_{i=1}^{m} \frac{1}{a_i}  
& \leq \sum_{i=1}^{m} \frac{1}{x_i}  - \left(  \sum_{i=m+1}^{n} \frac{1}{a_i} - \sum_{i=m+1}^{n} \frac{1}{x_i}  \right) \\
& < \sum_{i=1}^{m} \frac{1}{x_i}  \leq \sum_{i=1}^{n} \frac{1}{x_i} < \frac{p}{q}.
\end{align*} 
The induction hypothesis implies $x_i = a_i$ for all $i = 1,\ldots, m$, which is absurd.  
Thus,  $x_i = a_i$ for all $i = m+1,\ldots, n$, and 
\[
\sum_{i=1}^m \frac{1}{a_i} \leq \sum_{i=1}^m \frac{1}{x_i} < \sum_{i=1}^n \frac{1}{x_i} < \frac{p}{q}. 
\]   
The induction hypothesis again  implies $x_i = a_i$ for all $i = 1,\ldots, m$. 
This completes the proof.  
\end{proof}

\section{When is greedy best?}                 \label{Egyptian:section:cc}

It is a critical observation that the $n$-term greedy  underapproximation 
of a real number $\theta \in (0,1]$ 
is not always the unique best $n$-term Egyptian underapproximation, 
nor even a best  $n$-term Egyptian underapproximation.

Here are two examples for the case $n=2$.  The inequality
\[
\frac{1}{2} + \frac{1}{30} = \frac{8}{15} <  \frac{31}{58} = \frac{1}{2} + \frac{1}{29}
\] 
proves that  $(2,30)$ is the 2-term greedy underapproximation sequence for all $\theta$ 
in the interval 
\[
\frac{8}{15} < \theta \leq \frac{31}{58}.  
\] 
We prove (Theorem~\ref{Egyptian:theorem:underapproximation-a1=2}) 
that $(2,30)$ is a best 2-term greedy underapproximation sequence for all $\theta$ 
in this interval.  The equation  
\[
\frac{1}{2} + \frac{1}{30} = \frac{1}{3} + \frac{1}{5}  = \frac{8}{15} 
\]
shows that the best  2-term Egyptian underapproximation 
is not unique.

Similarly, the inequality
\[
 \frac{1}{3} + \frac{1}{17}   =  \frac{20}{51}  < \frac{19}{48} =  \frac{1}{3} + \frac{1}{16} 
\] 
proves that  $(3,17)$ is the 2-term greedy underapproximation sequence for all $\theta$ 
in the interval 
\[
 \frac{20}{51}  < \theta \leq \frac{19}{48} .  
\] 
The inequality 
\[
 \frac{1}{3} + \frac{1}{17}  < \frac{1}{4} + \frac{1}{7} =  \frac{11}{28}  < 
  \theta \leq \frac{19}{48} 
\]
proves that  $(3,17)$ is not a best 2-term Egyptian underapproximation 
of $\theta$ for all $\theta$ 
in the interval 
\[
 \frac{11}{28} < \theta \leq \frac{19}{48} .  
\] 
Theorem~\ref{Egyptian:theorem:underapproximation-a1=3} shows that $(4,7)$ 
is the best 2-term Egyptian underapproximation 
of $\theta$ for all $\theta$ in this interval.

\section{Best 2-term Egyptian underapproximations} 

In this section we describe best 2-term Egyptian underapproximations for 
$\theta \in (0,1]$.

For all integers $a_1 \geq 2$ we have the \emph{harmonic interval}\index{harmonic interval} 
\[
I(a_1) = \left(\frac{1}{a_1}, \frac{1}{a_1-1} \right] = \left(\frac{1}{a_1}, \frac{1}{a_1} + \frac{1}{a^2_1- a_1} \right]. 
\]
The intervals $I(a_1)$ are pairwise disjoint   and 
\[
(0,1] = \bigcup_{a_1=2}^{\infty} \left(\frac{1}{a_1}, \frac{1}{a_1-1} \right]. 
\]
For all integers $a_1 \geq 2$ and $a_2 \geq a_1^2-a_1+1$, 
we have the  \emph{harmonic subinterval}\index{harmonic subinterval}  
\[
J(a_1,a_2) = \left(\frac{1}{a_1} + \frac{1}{a_2}, \frac{1}{a_1} + \frac{1}{a_2-1} \right]. 
\]
By Corollary~\ref{Egyptian:corollary:GUAsequence-1}, the pair 
$(a_1,a_2)$ is the 2-term greedy underapproximation of $\theta$ for all $\theta \in J(a_1,a_2)$.

We have 
\[
J(a_1,a_2) \subseteq  \left(\frac{1}{a_1}, \frac{1}{a_1-1} \right]  = I(a_1).  
\]
The intervals $J(a_1,a_2)$  are pairwise disjoint.  
It follows from the identity  
\[
\left(\frac{1}{a_1} + \frac{1}{a_1^2 - a_1+1}, \frac{1}{a_1} + \frac{1}{a_1^2 - a_1} \right] 
= \left(\frac{1}{a_1} + \frac{1}{a_1^2 - a_1+1}, \frac{1}{a_1} + \frac{1}{a_1-1} \right]  
\]
that  
 \begin{align*}
I(a_1) &  = \bigcup_{a_2 = a_1^2-a_1+1}^{\infty} J(a_1,a_2) 
 = \bigcup_{a_2 = a_1^2-a_1+1}^{\infty} \left(\frac{1}{a_1} + \frac{1}{a_2}, \frac{1}{a_1} + \frac{1}{a_2-1} \right]. 
 \end{align*}
Thus, 
\[
(0,1] = \bigcup_{a_1=2}^{\infty} \quad  \bigcup_{a_2 = a_1^2-a_1+1}^{\infty} 
\left(\frac{1}{a_1} + \frac{1}{a_2}, \frac{1}{a_1} + \frac{1}{a_2-1} \right]. 
\]

The pair of integers $(x_1,x_2)$ with $2 \leq x_1 \leq x_2$ is not the 2-term greedy 
underapproximation sequence of some $\theta \in (0,1]$
if and only if $x_2 \leq x_1^2-x_1$.

The pair $(a_1,a_2)$ is not a best  2-term underapproximation sequence 
of  $\theta \in I(a_1,a_2)$ if and only if 
there exists a pair of positive integers $(x_1,x_2) $ with $x_1 \leq x_2$ such that 
\[
\frac{1}{a_1} + \frac{1}{a_2} < \frac{1}{x_1} + \frac{1}{x_2} <  \theta \leq \frac{1}{a_1} + \frac{1}{a_2-1}.
\]

The following Lemmata enable us to compute, for all integers $a_1 \geq 2$,  
the set of real numbers $\theta$  in the harmonic interval 
$I(a_1) = (1/a_1, 1/(a_1 -1)]$ for which the 2-term greedy underapproximation is not the 
unique best 2-term Egyptian underapproximation. 

\bl                                        \label{Egyptian:lemma:calculate-1} 
Let $a_1$ and $a_2$  be integers such that 
\[
a_1 \geq 2
\qqand  
a_2 \geq a_1(a_1-1)+1. 
\]
If $x_1$ and $x_2$ are  integers such that 
\[
2 \leq x_1 \leq x_2 \qqand (x_1,x_2) \neq (a_1, a_2) 
\] 
and 
\beq                                \label{Egyptian:estimate-x1-x2-0} 
\frac{1}{a_1} + \frac{1}{a_2} \leq \frac{1}{x_1} + \frac{1}{x_2} <  \frac{1}{a_1} + \frac{1}{a_2-1} 
\eeq
then 
\beq                                \label{Egyptian:estimate-x1-x2} 
a_1 + 1 \leq x_1 \leq 2a_1-1 \leq x_2 < \frac{a_1x_1}{x_1-a_1} 
\eeq
and
\beq                                  \label{Egyptian:greedy-a2x2-2} 
 x_2 \leq  a_2-1.
\eeq 
\el

\begin{proof}
We have    
\[
\frac{2}{x_2}  \leq \frac{1}{x_1} + \frac{1}{x_2} <  \frac{1}{a_1} + \frac{1}{a_2-1}  \leq \frac{1}{a_1-1}
\]
and so $2a_1 - 1 \leq x_2$.  
Similarly, 
\[
 \frac{1}{x_1} < \frac{1}{x_1} + \frac{1}{x_2} <    \frac{1}{a_1-1} 
\] 
implies $ a_1 \leq x_1$.  
If $a_1 = x_1$, then from~\eqref{Egyptian:estimate-x1-x2-0} we obtain  
\[
 \frac{1}{a_2} \leq  \frac{1}{x_2} <  \frac{1}{a_2-1} 
\]
and so $a_2 = x_2$, which contradicts $(x_1,x_2) \neq (a_1, a_2)$.   
It follows that  $a_1 +1 \leq x_1$.   

If $x_1 \geq 2a_1$, then 
\[
\frac{1}{a_1} + \frac{1}{a_2} \leq  \frac{1}{x_1} + \frac{1}{x_2} 
\leq \frac{2}{x_1} \leq \frac{1}{a_1} < \frac{1}{a_1} + \frac{1}{a_2}
\]
which is absurd. 
Therefore, 
\[
a_1 + 1 \leq x_1 \leq 2a_1-1 \leq x_2.
\]

The inequality 
\[
 \frac{1}{a_1}  < \frac{1}{a_1} + \frac{1}{a_2} \leq \frac{1}{x_1} + \frac{1}{x_2} 
 \]
implies 
\[
\frac{x_1-a_1}{a_1 x_1}  = \frac{1}{a_1}  -  \frac{1}{x_1} <  \frac{1}{x_2} 
\]
and so 
\[
2a_1-1 \leq x_2 < \frac{a_1 x_1} {x_1-a_1}. 
\]
This finishes the proof of~\eqref{Egyptian:estimate-x1-x2}.  

Finally, $a_1 < x_1$ implies 
\[
\frac{1}{a_1} + \frac{1}{a_2} \leq \frac{1}{x_1} + \frac{1}{x_2} < \frac{1}{a_1} + \frac{1}{x_2} 
\]
and so  $x_2 \leq a_2 -1$, which is~\eqref{Egyptian:greedy-a2x2-2}.  This completes the proof.  
\end{proof}

\bl                                                              \label{Egyptian:lemma:calculate-2}            
For all integers $a_1\geq 2$, there are $a_1-1$ integers $x_1$ that satisfy 
\[
a_1 + 1 \leq x_1 \leq 2a_1-1.
\]  
For each such integer $x_1$ there are 
\[
\left\lceil  \frac{a_1x_1}{x_1-a_1}\right\rceil  - 2a_1 + 2 \geq 1
\]
integers $x_2$ that satisfy 
\[
 2a_1-1 \leq x_2 < \frac{a_1x_1}{x_1-a_1}
\]
For all integers $a_1 \geq 2$, the set 
\beq                             \label{Egyptian:condition-2-term} 
X(a_1) = \left\{ (x_1,x_2) \in \N^2: a_1 + 1 \leq x_1 \leq 2a_1-1 \leq x_2 < \frac{a_1x_1}{x_1-a_1} \right\} 
\eeq 
is nonempty. 
\el

\begin{proof}
If $a_1\geq 2$, then $a_1+1 \leq 2a_1-1$.  There are $a_1-1 \geq 1$ integers $x_1$ such that 
$a_1+1 \leq x_1 \leq 2a_1-1$.  

If $x_1 \leq 2a_1-1$, then 
\[
(a_1-1)x_1 \leq (a_1 -1)(2a_1 -1) < a_1(2a_1-1).
\]
Equivalently, 
\[
(2a_1-1) (x_1  - a_1) <  a_1 x_1
\]
and so 
\[
 2a_1-1 < \frac{a_1x_1}{x_1-a_1}. 
\]
It follows that there are 
\[
\left\lceil  \frac{a_1x_1}{x_1-a_1}\right\rceil  - 2a_1 + 2 \geq 1
\]
integers $x_2$ such that 
\[
 2a_1-1 \leq x_2 < \frac{a_1x_1}{x_1-a_1}
\]
and so  the set $X(a_1)$ is nonempty.  
This completes the proof. 
\end{proof}

\bl                                                       \label{Egyptian:lemma:calculate-3} 
Let $a_1 \geq 2$.  
If $(x_1,x_2) \in X(a_1)$ and 
\beq                             \label{Egyptian:calculate-a2} 
a_2 =\left\lceil \left( \frac{1}{x_1} +  \frac{1}{x_2}  -  \frac{1}{a_1} \right)^{-1}  \right\rceil
\eeq
then
\beq                             \label{Egyptian:condition-2-term-delete} 
\frac{1}{a_1} + \frac{1}{a_2}  \leq \frac{1}{x_1}  + \frac{1}{x_2}  <  
 \frac{1}{a_1}  + \frac{1}{a_2 -1} 
\eeq
and the pairs  $(a_1,a_2)$ and $(x_1, x_2)$ are 2-term underapproximation 
sequences of $\theta$ for all 
\[
\theta \in \left(  \frac{1}{x_1}  + \frac{1}{x_2},  \frac{1}{a_1}  + \frac{1}{a_2 -1} \right].
\]
Moreover,  
\[
\frac{1}{a_1} + \frac{1}{a_2}=  \frac{1}{x_1} + \frac{1}{x_2} 
\]
if and only if 
\[
a_2 =  \left(    \frac{1}{x_1} +  \frac{1}{x_2}  -  \frac{1}{a_1} \right)^{-1}.  
\]
\el

\begin{proof} 
If $(x_1,x_2) \in X(a_1)$, then  
\[
a_2 -1 <  \left( \frac{1}{x_1} +  \frac{1}{x_2}  -  \frac{1}{a_1} \right)^{-1}  \leq a_2 
\]
and 
\[
  \frac{1}{a_1} + \frac{1}{a_2} \leq \frac{1}{x_1} +  \frac{1}{x_2}   <    \frac{1}{a_1} + \frac{1}{a_2-1}. 
\]
This proves~\eqref{Egyptian:condition-2-term-delete}.  
The remaining statements are immediate consequences.  
\end{proof}

It is important to note that the integer $a_2$ computed from~\eqref{Egyptian:calculate-a2}  
does not necessarily satisfy the inequality $a_2 \geq a_1^2 - a_1 + 1$.  Thus,  
$(x_1,x_2)$ is an equal or better 2-term underapproximation than $(a_1,a_2)$ 
for all $\theta > 1/x_1 + 1/x_2$, but 
$(a_1,a_2)$ is not necessarily a 2-term greedy underapproximation.

Here are three examples in the case $a_1 = 5$.  We have $a_1^2-a_1+1=21$ and 
\[
 \left( \frac{1}{5} ,   \frac{1}{4}  \right] = I(5) = \bigcup_{a_2 = 21}^{\infty} J(5,a_2) 
= \bigcup_{a_2 = 21}^{\infty} \left( \frac{1}{5} + \frac{1}{a_2},   \frac{1}{5} + \frac{1}{a_2-1} \right]. 
\]
For all $a_2 \geq 21$, the pair $(a_1,a_2) = (5,a_2)$ is the 2-term greedy underapproximation  
of $\theta$ of all $\theta$ is in the harmonic subinterval $J(5,a_2)$. 
From~\eqref{Egyptian:condition-2-term} we obtain the inequality that determines the set  $X(5)$: 
\[
6 \leq x_1 \leq 9 \leq x_2 < \frac{5x_1}{x_1-5}.
\]
The set $X(5)$ contains the pairs  $(x_1,x_2) = (7,10)$,  $(9,11)$, and $(6,9)$.

The pair $(7,10) \in X(5)$ generates the integer  
\[
a_2 = 24 >   \left( \frac{1}{7} +  \frac{1}{10}  -  \frac{1}{5} \right)^{-1} = \frac{70}{3} > 23
\]
and  $24 = a_2  \geq   21$.    
The pair $(5,24)$ is the 2-term greedy underapproximation sequence of all $\theta \in J(5,24)$.  
We have 
\[
  \frac{29}{120} =  \frac{1}{5} + \frac{1}{24} < \frac{1}{7} +  \frac{1}{10}  = \frac{17}{70} 
  <    \frac{1}{5} + \frac{1}{23} =  \frac{28}{115}. 
\]
Thus, the pair $(7,10)$ is a better 2-term underapproximation sequence of $\theta$  than 
the 2-term greedy underapproximation sequence $(5,24)$  for all  
\[
\theta \in \left(  \frac{17}{70},  \frac{28}{115} \right] \subseteq \left(  \frac{29}{120},  \frac{28}{115} \right] 
= J(5,24). 
\]

The pair $(9,11) \in X(5)$ generates the integer  
\[
a_2 =  495 =   \left( \frac{1}{9} +  \frac{1}{11}  -  \frac{1}{5} \right)^{-1} 
\]
and $495 = a_2  \geq   21$.  
The pair $(5,495)$ is the 2-term greedy underapproximation sequence 
for all $\theta \in J(5,495)$.   
For all  
\[
\theta \in \left( \frac{20}{99},\frac{499}{2470}  \right] 
\]
we have 
\[
 \frac{1}{5} + \frac{1}{495} = \frac{1}{9} +  \frac{1}{11}  = \frac{20}{99} 
  <   \theta \leq \frac{499}{2470} =   \frac{1}{5} + \frac{1}{494}.
\]
and the pairs $(5,495)$ and$(9, 11)$ give equal 2-term underapproximations.

The pair $ (6,9) \in X(5)$ generates the integer    
\[
a_2 = 13 =   \left\lceil    \frac{90}{7}  \right\rceil  
= \left\lceil   \left( \frac{1}{6} +  \frac{1}{9}  -  \frac{1}{5} \right)^{-1} \right\rceil >  \frac{90}{7}.
\]
However, $a_2 = 13 < 21$ and  $(5,13)$ is not a 2-term greedy underapproximation sequence.

\section{Best 2-term underapproximations for $a_1=2$ and $a_1 = 3$}           \label{Egyptian:section:dd}
In this section we compute all real numbers $\theta$ in the harmonic intervals 
$I(2)$ and  $I(3)$  whose 2-term greedy approximation sequences do not give best approximations 
or unique best approximations.  

\bt               \label{Egyptian:theorem:underapproximation-a1=2} 
Let $a_1=2$ and $a_2 \geq 3$.  
The 2-term greedy underapproximation sequence $(2,a_2)$ is a best 
2-term Egyptian underapproximation sequence of $\theta$ 
for all $\theta$ in the harmonic subinterval 
\[
J(2,a_2) =  \left( \frac{1}{2} + \frac{1}{a_2},  \frac{1}{2} + \frac{1}{a_2-1} \right]. 
\]   
Consider the harmonic subintervals 
\begin{align*}
J(2,6) & =  \left( \frac{2}{3},\frac{7}{10} \right] , \qquad 
J(2,12)  =  \left( \frac{7}{12},\frac{13}{22} \right], \qquad 
J(2,30)   =  \left( \frac{8}{15},\frac{31}{58} \right].   
\end{align*}  
\benum
\item[(i)]
For all $\theta \in J(2,6)$, the pairs $(2,6)$ and   $(3,3)$ are best 2-term  underapproximations of $\theta$, 
and are the only  best 2-term  underapproximations of $\theta$.  
\item[(ii)]
For all $\theta \in J(2,12)  $,  the pairs $(2,12)$ and   $(3,4)$ are best 2-term  underapproximations of $\theta$, 
and are the only  best 2-term  underapproximations of $\theta$.   

\item[(iii)]
For all $\theta \in J(2,30)  $,  the pairs $(2,30)$ and   $(3,5)$ are best 2-term  underapproximations of $\theta$, 
and are the only  best 2-term  underapproximations of $\theta$.   
\item[(iv)]
For all $\theta \in I(2) = (1/2,1]$ such that $\theta \notin    J(2,6)   \cup J(2,12)   \cup J(2,30) $, 
 the pair $ (a_1,a_2)$  is the unique best 2-term  underapproximation of $\theta$. 
 \eenum
\et

\begin{proof}
If  $a_1=2$, then inequality~\eqref{Egyptian:estimate-x1-x2} is simply 
\[
3 = x_1  \leq x_2 < 6
\]
and so $x_2 = 3, 4$, or 5.  
If $x_2 = 3$, then 
\[
a_2 = \left( \frac{1}{3} + \frac{1}{3} - \frac{1}{2} \right)^{-1} = 6 
\]
and 
\[
 \frac{1}{2} +  \frac{1}{6} = \frac{1}{3} + \frac{1}{3} = \frac{2}{3}. 
\]
If $x_2 = 4$, then 
\[
a_2 = \left( \frac{1}{3} + \frac{1}{4} - \frac{1}{2} \right)^{-1} = 12 
\]
and 
\[
 \frac{1}{2} +  \frac{1}{12} = \frac{1}{3} + \frac{1}{4} = \frac{7}{12}. 
\] 
If $x_2 = 5$, then 
\[
a_2 = \left( \frac{1}{3} + \frac{1}{5} - \frac{1}{2} \right)^{-1} = 30 
\]
and 
\[
 \frac{1}{2} +  \frac{1}{30} = \frac{1}{3} + \frac{1}{5} = \frac{8}{15}. 
\]
The only solutions $(x_1,x_2) \neq (2,a_2)$ of the diophantine inequality  
\beq                  \label{Egyptian:2term-ineq}
\frac{1}{2} + \frac{1}{a_2}  \leq \frac{1}{x_1}  + \frac{1}{x_2} < \frac{1}{2} + \frac{1}{a_2 -1} 
\eeq
are $(x_1,x_2) = (3,3)$, $(3,4)$, and $(3,5)$. 
This completes the proof.   
\end{proof}

\bt               \label{Egyptian:theorem:underapproximation-a1=3} 
Let $a_1=3$ and let $\theta \in I(3) = (1/3,1/2]$.  
The 2-term greedy underapproximation of $\theta$ is   a 
best 2-term  Egyptian underapproximation if and only if 
\[
\theta \notin  \left( \frac{9}{20} , \frac{11}{24} \right]  \cup \left( \frac{11}{28} , \frac{19}{48} \right]. 
\]
The 2-term greedy underapproximation of $\theta$ is   a 
best 2-term  Egyptian underapproximation  but not 
the unique best 2-term  Egyptian underapproximation if and only if 
$\theta \in J(3,a_2)$ for 
\[
a_2 \in \{12, 15, 24, 30, 36, 60, 105, 132 \}. 
\]
\et

\begin{proof}
For $a_1=3$, inequality~\eqref{Egyptian:estimate-x1-x2}  gives 
\[
4 \leq x_1 \leq 5 \leq x_2  < \frac{3x_1}{x_1-3}.
\]
Thus, a complete list of the 10 solutions $(3,a_2) \neq (x_1,x_2)$ of the diophantine inequality  
\beq                  \label{Egyptian:underapproximation-3-equal}
\frac{1}{3} + \frac{1}{a_2}  \leq \frac{1}{x_1}  + \frac{1}{x_2} < \frac{1}{3} + \frac{1}{a_2 -1} 
\eeq  
is the following:  
\beq                                    \label{Egyptian:underapproximation-3-table}
\begin{tabular}{|r | r | r  | r  | r  | r  | r  | r  | r  | r  | r |} \hline 
$x_1$ & 4 & 4  & 4  & 4  & 4  & 4  & 4  & 5  & 5  & 5 \\ 
$x_2$ & 5 & 6 & 7 & 8 & 9 & 10 & 11 & 5 & 6 & 7 \\  
$a_2$ & 9 & 12 & 17 & 24 & 36 & 60 & 132 & 15  & 30  &  105  \\ \hline
\end{tabular}
\eeq
We have strict inequality 
\[
\frac{1}{3} + \frac{1}{a_2}  < \frac{1}{x_1}  + \frac{1}{x_2} < \frac{1}{3} + \frac{1}{a_2 -1} 
\]
only if either $a_2 = 9 $ and $(x_1,x_2) = (4,5)$ or $a_2 = 17$ and $(x_1,x_2) = (4,7)$.  
Note that 
\[
\frac{1}{4} + \frac{1}{5}  = \frac{9}{20} \qqand \frac{1}{4}  + \frac{1}{7} = \frac{11}{28}. 
\]
The pair $(4,5)$ is the unique best 2-term  underapproximation of all $\theta$ such that 
\[
\theta \in  \left( \frac{9}{20} , \frac{11}{24} \right]  \subseteq  \left( \frac{4}{9}, \frac{11}{24} \right]  = J(3,9).  
 \]
The pair $(4,7)$ is the unique best 2-term  underapproximation of all $\theta$ such that 
\[
\theta \in \left( \frac{11}{28} , \frac{19}{48} \right] \subseteq \left( \frac{20}{51} , \frac{19}{48} \right] = J(3,17).
\]

The 8 solutions $(3,a_2) \neq (x_1,x_2)$ with $a_2 \geq 7$ and $4 \leq x_1 \leq x_2$ 
of the diophantine equation  
\[
\frac{1}{3} + \frac{1}{a_2}  = \frac{1}{x_1}  + \frac{1}{x_2}  
\]
are
\[
\frac{1}{3}  +  \frac{1}{ 12} = \frac{1}{4}  +  \frac{1}{ 6} =\frac{5}{12 } 
\]
\[
\frac{1}{3}  +  \frac{1}{15 } = \frac{1}{5}  +  \frac{1}{5} =  \frac{2}{5} 
\]
\[
\frac{1}{3}  +  \frac{1}{24 } = \frac{1}{4}  +  \frac{1}{ 8} =  \frac{3}{8 } 
\]
\[
\frac{1}{3}  +  \frac{1}{30} = \frac{1}{5}  +  \frac{1}{ 6} =\frac{11}{30} 
\]
\[
\frac{1}{3}  +  \frac{1}{36 } = \frac{1}{4}  +  \frac{1}{9 } =  \frac{13}{ 36} 
\]
\[
\frac{1}{3}  +  \frac{1}{60 } = \frac{1}{4}  +  \frac{1}{10 } =  \frac{7}{ 20} 
\]
\[
\frac{1}{3}  +  \frac{1}{105} = \frac{1}{5}  +  \frac{1}{7} =  \frac{12}{35 }. 
\]
\[
\frac{1}{3}  +  \frac{1}{ 132} = \frac{1}{4}  +  \frac{1}{11 } =  \frac{15}{44 }. 
\]

This completes the proof. 
\end{proof}

\section{Open Problems}
\benum

\item
Consider real numbers $\theta \in (0,1]$ whose infinite greedy underapproximation sequence 
$(a_i)_{i=1}^{\infty}$ has the property that $(a_i)_{i=1}^n$ is the unique best underapproximation 
of $\theta$ for all positive integers $n$. By Theorem~\ref{Egyptian:theorem:pq-greedy}, 
every rational number of the form 
$p/q$ where $p$ divides $q+1$ has this property.  
Do other rational numbers have this property?   
Do there exist irrational numbers with this property?

\item 
Let $\theta \in (0,1]$, let $n \geq 3$, and let $(a_i)_{i=1}^n \in U_n(\theta)$ 
be the $n$-term  greedy underapproximation sequence of $\theta$. 
\benum
\item
Do  there exist sequences $(x_i)_{i=1}^n \in U_n(\theta)$ such that 
$(a_i)_{i=1}^n \neq (x_i)_{i=1}^n$ and 
\[
 \sum_{i=1}^n \frac{1}{a_i} <  \sum_{i=1}^n \frac{1}{x_i} < \theta?
\]
How many such sequences are there? 
\item
Do  there exist sequences $(x_1,\ldots, x_n) \in U_n(\theta)$ such that 
$(a_1,\ldots, a_n)  \neq (x_1,\ldots, x_n)$ and 
\[
 \sum_{i=1}^n \frac{1}{a_i} =  \sum_{i=1}^n \frac{1}{x_i} < \theta?
\]
How many such sequences are there?  
\item
Can we identify and understand counterexamples to unique best $n$-term 
underapproximation by the greedy algorithm?
\eenum

\item 
Let $n \geq 3$.  
Is there an efficient algorithm to compute the best $n$-term underapproximation 
sequence of a real number $\theta \in (0,1]$?

\item
Let $\theta \in (0,1]$. 
Erd\H os and Graham~\cite[p.31]{erdo-grah80} asserted (without proof or reference to any publication)  
that for every rational number $\theta$ 
there exists an integer $n_0 = n_0(\theta)$ such that, 
for all $n \geq n_0 + 1$,  
\[
u_n(\theta) = u_{n_0}(\theta) + u_{n-n_0}\left( \theta - u_{n_0}(\theta) \right)
\]
and the best $(n-n_0)$-term underapproximation $u_{n-n_0}\left( \theta - u_{n_0}(\theta) \right)$ 
is always constructed by the greedy algorithm.  They also wrote,  ``It is not difficult to construct 
irrationals for which the result fails.''
Prove or disprove these statements.

\item
Let $A$ be a nonempty set of positive integers and let 
\[
\frac{1}{A} = \left\{ \frac{1}{x} : x \in A \right\} 
\] 
be the set of Egyptian fractions with denominators in $A$.  
An \emph{$n$-term $A$-underapproximation}\index{Egyptian underapproximation!$A$} 
of  $\theta$ is a sum of $n$ not necessarily distinct Egyptian fractions in  $1/A$ 
that is strictly less than $\theta$.  
Let
\[
u_{n,A}(\theta) = \sup\left\{ \sum_{i=1}^n \frac{1}{x_i}: (x_1,\ldots, x_n) \in A^n, x_1 \geq \cdots \geq x_n, 
\sum_{i=1}^n \frac{1}{x_i} < \theta \right\}. 
\]
An $n$-term $A$-underapproximation $ \sum_{i=1}^n \frac{1}{x_i} < \theta$ is 
\emph{best}\index{Egyptian underapproximation!best} if  
\[
u_n(\theta) =  \sum_{i=1}^n \frac{1}{x_i}.
\]
For what real numbers $\theta \in (0,1)$ does the greedy algorithm restricted to $A$ 
give a best $n$-term underapproximation?
\eenum

\appendix                                                 \label{Egyptian:appendix}
\section{Proof of an inequality}

The following proof is due to Ambro and Barc\u{a}u~\cite{ambr-barc15}.
\bt                       
Let $(u_i)_{i=1}^n$ and $(v_i)_{i=1}^n$ be distinct sequences of positive numbers such that      
\beq                                    \label{Muirhead:Muirhead-AddMult-Ineq0}
\prod_{i=1}^k v_i \leq \prod_{i=1}^k u_i 
\eeq
for all $k = 1,\ldots, n$.   If  $(u_i)_{i=1}^n$ and $(v_i)_{i=1}^n$ are decreasing, then 
\beq                                \label{Muirhead:Muirhead-AddMult-Ineq1}
\sum_{i=1}^n v_i < \sum_{i=1}^n u_i.
\eeq
If  $(u_i)_{i=1}^n$ and $(v_i)_{i=1}^n$ are increasing, then 
\beq                                \label{Muirhead:Muirhead-AddMult-Ineq2}
\sum_{i=1}^n \frac{1}{u_i}  < \sum_{i=1}^n  \frac{1}{v_i} .
\eeq
\et

\begin{proof}
Let $(u_i)_{i=1}^n$ and $(v_i)_{i=1}^n$ be decreasing sequences that 
satisfy the product inequality~\eqref{Muirhead:Muirhead-AddMult-Ineq0}.  
The proof of inequality~\eqref{Muirhead:Muirhead-AddMult-Ineq1} is by induction on $n$.  
The case $n=1$ is simply the assertion 
that if $v_1 \leq u_1$ and $v_1 \neq u_1$, then $v_1 < u_1$.  

Let $n \geq 2$ and assume  inequality~\eqref{Muirhead:Muirhead-AddMult-Ineq1} 
is true for sequences of length less than $n$.  If   $u_i = v_i$ for some $i$, the inequality follows 
from the case $n-1$.  If $v_i < u_i$ for all $i = 1,\ldots, n$, then $\sum_{i=1}^n v_i < \sum_{i=1}^n u_i$.  
Thus, we can assume that $u_i \neq v_i$ for all  $i$ 
and that $u_j < v_j$ for some $j$.  

Let $\ell$ be the smallest $j$ such that $u_j < v_j$. 
Inequality~\eqref{Muirhead:Muirhead-AddMult-Ineq0} with $k=1$ gives $v_1 \leq u_1$ 
and so $v_1 < u_1$.  Therefore, $\ell \geq 2$ and 
\[
u_{\ell} <  v_{\ell} \leq  v_{\ell-1} < u_{\ell -1}.
\]
Let   
\beq           \label{Egyptian:Muir-a}
t = \min \left( \frac{u_{\ell-1}}{v_{\ell-1}}, \frac{ v_{\ell} }{ u_{\ell}} \right). 
\eeq
We have $ v_{\ell}  \leq v_{\ell-1}$ and 
\beq           \label{Egyptian:Muir-b}
1 < t < t^2 \leq  \left( \frac{u_{\ell-1}}{v_{\ell-1}} \right) \left(  \frac{ v_{\ell} }{ u_{\ell}} \right) 
= \left( \frac{u_{\ell-1}} { u_{\ell}}  \right) \left(  \frac{ v_{\ell} }  {v_{\ell-1}}\right) \leq  \frac{u_{\ell-1}} { u_{\ell}}. 
\eeq

Define the sequence $(u'_i)_{i = 1}^n$ as follows:
\begin{align*}
u'_i & = u_i  \quad \text{ if $i \neq \ell-1,\ell$} \\
u'_{\ell -1} & = \frac{u_{\ell-1}} {t} \\
u'_{\ell} &= t u_{\ell} . 
\end{align*}
Inequalities~\eqref{Egyptian:Muir-a} and~\eqref{Egyptian:Muir-b} imply that 
\[
 u_{\ell-1} >   \frac{u_{\ell-1}}{t} = u'_{\ell-1}   \geq  u'_{\ell}  =  tu_{\ell}   > u_{\ell} 
\]
and so 
\[
u'_1 \geq \cdots \geq u'_{\ell-2} \geq u_{\ell-1} >   u'_{\ell-1} 
 \geq u'_{\ell}  > u_{\ell} > u'_{\ell+1} \geq \cdots \geq u'_n. 
\]
Thus, the sequence $(u'_i)_{i=1}^n$ is decreasing.  

We shall prove product inequality~\eqref{Muirhead:Muirhead-AddMult-Ineq0} 
for the  decreasing sequences  $(u'_i)_{i=1}^n$ and $(v_i)_{i=1}^n$. 
For $k = 1,\ldots, \ell-2$ we have 
\[
 \prod_{i=1}^k v_i \leq \prod_{i=1}^k u_i =  \prod_{i=1}^k u'_i. 
\]
For $k = \ell  ,\ldots, n$, the identity  
\[
u'_{\ell-1} u'_{\ell} = \left(  \frac{u_{\ell-1}} {t} \right) \left(  t u_{\ell} \right) = u_{\ell-1} u_{\ell} 
\]
 implies  
\[
\prod_{i=1}^k v_i \leq   \prod_{i=1}^k u_i 
  = \left( \prod_{i=1}^{\ell-2}  u_i  \right) \left( u'_{\ell-1} u'_{\ell} \right) 
\left( \prod_{i= \ell+1}^k u_i  \right)  = \prod_{i=1}^k u'_i.  
\]
Let $k = \ell -1$.  We have $t \leq  u_{\ell-1}/v_{\ell-1}$ and so 
$ v_{\ell-1} \leq  u_{\ell-1}/t = u'_{\ell-1}$.   It follows that 
\[
 \prod_{i=1}^{\ell -1}  v_i  = v_{\ell-1}  \prod_{i=1}^{\ell -2}  v_i \leq  u'_{\ell-1}  \prod_{i=1}^{\ell -2} u'_i 
 =  \prod_{i=1}^{\ell -1} u'_i. 
\]
This proves~\eqref{Muirhead:Muirhead-AddMult-Ineq0} 
for the sequences  $(u'_i)_{i=1}^n$ and $(v_i)_{i=1}^n$. 
If $t = u_{\ell-1}/v_{\ell-1}$, then $ v_{\ell-1}  = u'_{\ell-1}$. 
If $t =  v_{\ell} / u_{\ell}$, then $ v_{\ell}  = u'_{\ell}$. 
In both cases, the induction hypothesis implies 
\[
\sum_{i=1}^n v_i < \sum_{i=1}^n u'_i. 
\]

We have $t <  u_{\ell-1}/u_{\ell}$ from inequality~\eqref{Egyptian:Muir-b} and so   
\[
(t-1)u_{\ell} < \left(1-\frac{1}{t}\right) u_{\ell-1}. 
\]
Equivalently, 
\[
u'_{\ell-1} +  u'_{\ell} = \frac{u_{\ell-1}}{t} + t u_{\ell} < u_{\ell-1} +  u_{\ell} 
\]
and
\begin{align*}
\sum_{i=1}^n u'_i 
& = \left(  \sum_{\substack{i=1 \\ i \neq \ell-1, \ell}}^n u'_i \right) + u'_{\ell-1} +  u'_{\ell}  \\
& <  \left(  \sum_{\substack{i=1 \\ i \neq \ell-1, \ell}}^n u'_i \right)  + u_{\ell-1} +  u_{\ell}  =  \sum_{i=1}^n u_i.
\end{align*} 
This proves the Theorem for decreasing sequences.

If $(u_i)_{i=1}^n$ and $(v_i)_{i=1}^n$ are increasing sequences 
that satisfy~\eqref{Muirhead:Muirhead-AddMult-Ineq0}, then  
 $(1/u_i)_{i=1}^n$ and $(1/v_i)_{i=1}^n$ are decreasing sequences 
 such that 
 \[
 \prod_{i=1}^k \frac{1}{u_i} \leq \prod_{i=1}^k \frac{1}{v_i}
\]
for all $k = 1,\ldots, n$, and we obtain inequality~\eqref{Muirhead:Muirhead-AddMult-Ineq2}.
This completes the proof.  
\end{proof}

Note that inequality~\eqref{Muirhead:Muirhead-AddMult-Ineq0}
does not imply inequality~\eqref{Muirhead:Muirhead-AddMult-Ineq1} 
if the sequences $(u_i)_{i=1}^n$ and $(v_i)_{i=1}^n$ are increasing. 
For example, the increasing sequences  $(v_1,v_2) = (1,7)$ and $(u_1,u_2) = ( 2,4) $ satisfy 
\[
1 = v_1 < u_1 = 2 \qqand  7 = v_1  v_2 < u_1  u_2 = 8 
\]
but 
\[
8 =  v_1 + v_2 > u_1 + u_2 = 6.
\]

\def\cprime{$'$} \def\cprime{$'$} \def\cprime{$'$}
\providecommand{\bysame}{\leavevmode\hbox to3em{\hrulefill}\thinspace}
\providecommand{\MR}{\relax\ifhmode\unskip\space\fi MR }
\providecommand{\MRhref}[2]{
  \href{http://www.ams.org/mathscinet-getitem?mr=#1}{#2}
}
\providecommand{\href}[2]{#2}

\end{document}